\documentclass[conference]{IEEEtran}
\IEEEoverridecommandlockouts
\usepackage{cite}

\ifCLASSINFOpdf
  \usepackage[pdftex]{graphicx}
\else
\fi
\usepackage{amsmath,amssymb,amsfonts}
\usepackage{algorithmic}
\usepackage{xcolor}

\usepackage{multicol}

\begin{document}
%
\title{Cloud-AC-OPF: Model Reduction Technique for Multi-Scenario Optimal Power Flow via Chance-Constrained Optimization}






%
\author{
\IEEEauthorblockN{Vladimir Frolov\IEEEauthorrefmark{1}, Line Roald\IEEEauthorrefmark{2}  and Michael Chertkov\IEEEauthorrefmark{1}\IEEEauthorrefmark{3}\IEEEauthorrefmark{4}}

\IEEEauthorblockA{\IEEEauthorrefmark{1}Skolkovo Institute of Science and Technology, Nobel Street 3, Moscow Region, 143026, Russia}
\IEEEauthorblockA{\IEEEauthorrefmark{2}Electrical \& Computer Engineering, University of Wisconsin, 1415 Engineering Drive, Madison, WI 53706}
\IEEEauthorblockA{\IEEEauthorrefmark{3}Center for Nonlinear Studies and Theoretical Division, T-4, LANL, Los Alamos, NM 87545}
\IEEEauthorblockA{\IEEEauthorrefmark{4}Program in Applied Mathematics, University of Arizona, Tucson, AZ 85721}

\thanks{The work at LANL was carried out under the auspices of the National Nuclear Security Administration of the U.S. Department of Energy under Contract No. DE-AC52-06NA25396. The work was partially supported by DOE/OE/GMLC and LANL/LDRD/CNLS projects}
}


\maketitle
\begin{abstract}
Many practical planning and operational applications in power systems require simultaneous consideration of a large number of operating conditions or Multi-Scenario AC-Optimal Power Flow (MS-AC-OPF) solution. However, when the number of exogenously prescribed conditions is large, solving the problem as a collection of single AC-OPFs may be time-consuming or simply intractable computationally. In this paper, we suggest a model reduction approach, coined Cloud-AC-OPF, which replaces a collection of samples by their compact representation in terms of mean and standard deviation. Instead of determining an optimal generation dispatch for each sample individually, we parametrize the generation dispatch as an affine function. The Cloud-AC-OPF is mathematically similar to a generalized Chance-Constrained AC-OPF (CC-AC-OPF) of the type recently discussed in the literature, but conceptually different as it discusses applications to long-term planning. We further propose a tractable formulation and implementation, and illustrate our construction on the example of 30-bus IEEE model.
\end{abstract}

%
\IEEEpeerreviewmaketitle

\section{Introduction}

A variety of problems in power systems require consideration of a large number of operating conditions, corresponding to uncertain or time-varying renewable generation and load, different weather conditions or different economic situations. Each operating condition is usually well represented by a combination of the AC power flow equations to model the power flow physics, decision variables to model controllable variables such as the generation dispatch, and a set of parameters to model a particular operating condition such as realization of the load or renewable energy generation. Each operating condition hence gives rise to a unique instance of the AC Optimal Power Flow (AC-OPF) problem, and considering all operating conditions together results in what we will refer to as the Multi-Scenario AC-OPF (MS-AC-OPF).


While significant efforts have been invested into designing and implementing solution methods \cite{EES-012} 
and solvers \cite{MATPOWER, 06WB} for the AC-OPF problem, solving the MS-AC-OPF can be a challenging and time-consuming task. Even if the operating conditions, and thus decisions, are separable so that the optimization results in solving collection of single AC-OPF problems (one per-sample). For example, representing a year by 35,000 independent scenarios, each corresponding to a 15 min interval (with no uncertainty), will take about $10^5$ sec which is more than a day of computing (assuming that each problem is solved in 3 sec). An even more relevant case is the situation of a two-stage program, where we make decisions on the first stage (such as investment decisions) such that many operational situations contributing the MS-AC-OPF on the second stage are \emph{simultaneously} feasible \cite{93GCCP,14FBCa,14FBCb,2015MW}. In particular, as penetration of the renewable technology grows, also resulting in increase of uncertainty in generation and power flows, it becomes important to consider carefully how to identify the scenarios to represent and develop tractable methods to solve the resulting problems. 

In this manuscript, we suggest an alternative approach. Instead of considering a large number of samples to represent time-varying loads and renewable generation, we aggregate the samples into so-called \emph{scenario clouds}. For each cloud, we define a mean and covariance for the time-varying parameters. We then parametrize the generation dispatch as an affine function of the random parameters, which is a conservative (sub-optimal), but feasible choice for generation control. We use this to formulate what we refer to as the Cloud-AC-OPF, which is a mathematically similar, but conceptually different variant of a chance-constrained AC-OPF (CC-AC-OPF) \cite{13ROKA,14BCH,18RA,venzke2018chance}
with a prescribed robustness level. Hence, if we establish that the Cloud-AC-OPF is a good approximation of the MS-AC-OPF, we get algorithms which are capable of solving efficiently complex problems where the MS-AC-OPF represent the second (inner) stage, such as the decomposition methods developed to handle contingencies in the chance-constrained unit commitment in \cite{sundar2016unit}.

Contributions of this paper is three-fold. First, we propose the model reduction of the full MS-AC-OPF to the Cloud-AC-OPF. Second, we provide an efficient solution algorithm to solve the Cloud-AC-OPF. Third, we investigate how well the Cloud-AC-OPF approximates the true cost of the MS-AC-OPF. To keep our analysis clear and concise, we consider in this manuscript the case with limited variability corresponding to a single cloud and test the Cloud-AC-OPF performance against the MS-AC-OPF solved independently for each scenario. The analysis includes different parametrizations of the generation dispatch and different system loading levels. Future work will consider a larger range of variability (representative of scenarios occuring over a longer time horizon), which will lead to the definition of multiple scenario clouds with corresponding Cloud-AC-OPFs (Multi-Cloud-AC-OPF generalization).

\section{Why and How of the Model Reduction}
In this manuscript we focus on a problem which represent one piece of a bigger problem that we aim to address in the future -- namely long-term investment planning in transmission level of power system taking into account uncertainty of future operational conditions. The main difficulty in addressing this bigger problem is related to the significant variability and uncertainty of future operating conditions. One approach for taking this variability into account consists in modeling all possible operational conditions and incorporating these into the planning problem. This results in a two-stage optimization problem where the first stage optimizes over capacity/investment decisions, and the second stage optimizes the operational decisions (one set per scenario) subject to constraints given by the first stage \cite{14FBCa,14FBCb,16FPBBC,17FC}.

\subsection{Multi-Scenario AC-OPF}
In this manuscript, we consider the second stage of this future formulation. A typical formulation of this problem is the MS-AC-OPF, which can be stated (schematically) as follows:

\begin{subequations}
\begin{eqnarray}
\mbox{MS-AC-OPF}(x_u^{(a)}|\forall a)&=&\min_{x_c^{(a)}} \mbox{Cost}(x_c^{(a)})\label{eq:MS-AC-obj}\\
\mbox{s.t. }\forall a: & & \mbox{AC-PF}(x_c^{(a)},x_u^{(a)})=0\label{eq:AC-PF}\\ && \mbox{Constr}(x_c^{(a)},x_u^{(a)})\leq 0\label{eq:AC-ineq}
\end{eqnarray}
\label{eq:MS-AC-OPF}
\end{subequations}
Here, $x_c^{(a)},x_u^{(a)}$ are controlled and, respectively, uncontrolled state variables.
Values of the traditional generator dispatches are standard controlled variables. The output of renewable generators as well as consumption of many (aggregated) loads are examples of the uncontrolled variables, that are given input parameters to the problem. Some state variables, such as voltages and phases at all buses of the system (except for the slack bus) can be considered controlled, but are uniquely determined by the AC-PF equations  defined by (\ref{eq:AC-PF}). 
In (\ref{eq:MS-AC-OPF}), $a$ indexes a sample, $a=1,\cdots,M$; the objective \eqref{eq:MS-AC-obj} accounts for the cost of the traditional generation; the inequalities in the conditions of (\ref{eq:AC-ineq}) express the constraints on line flows, voltages, etc, which are introduced to enforce safe operations. 

\begin{figure}[h!]
\centering
\includegraphics[width=0.45\textwidth]{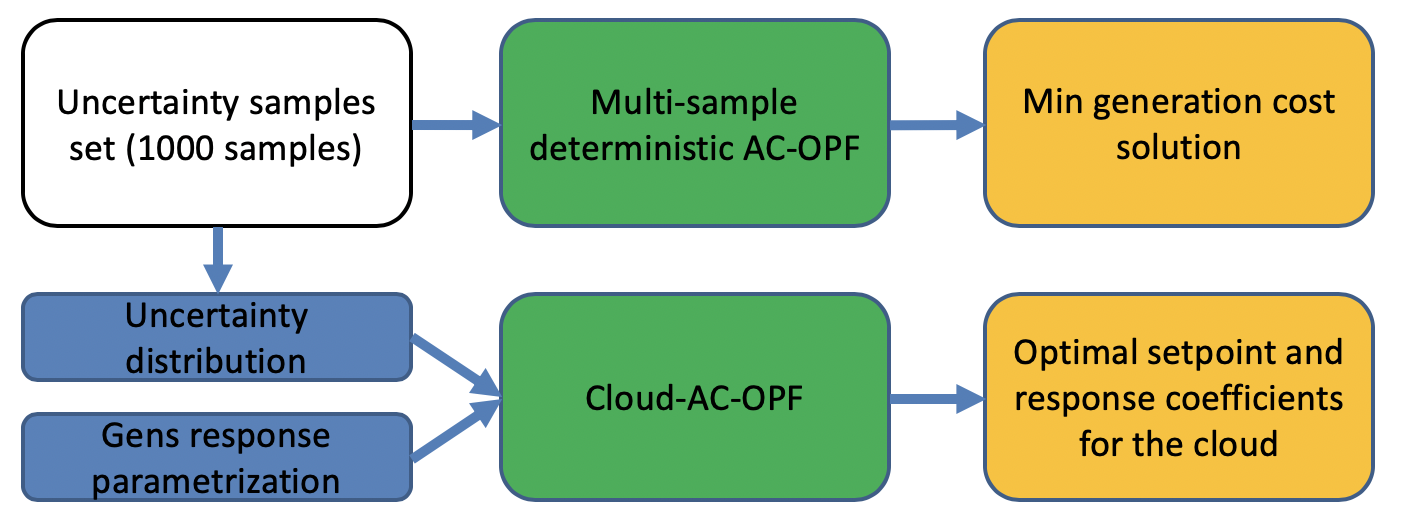}
\caption{Computational complexity reduction scheme}
\label{reduction_schematic}
\end{figure}

In many cases of practical interest the number of samples, $M$, can be too large to allow for sufficient accuracy (when approached in a brute-force fashion). A way to bypass this difficulty is in using stochastic methodology to describe the uncertainty set and also in representing system response to the uncertainty in a reduced, low-parametric way.  High-level logic of the reduction scheme employed in this manuscript is illustrated in  Fig.{\ref{reduction_schematic}}. 



\subsection{Uncertainty Modelling}
\label{subsec:uncertainty_modeling}
The key idea of our model reduction approach is to represent a group of similar samples by a single cloud. For example, a cloud can be formed by summer off-peak or peak operational conditions on hourly basis. Time correlation between uncertain variables is not considered. Mathematically cloud can be defined using clustering of historical multidimensional state vectors on the given time basis.
Then each sample $x_u^{(a)}$ represents an uncontrolled configuration assigned to the cloud. Given all samples in the cloud $x_u^{(a)}$ for $a=1,...,M$, we can calculate the mean $\bar{x}_u$ and the standard deviation $\sigma_u$ of the cloud for each uncertain bus according to:
\begin{align}
\bar{x}_u=1/M\sum_{a=1..M}x_u^{(a)},\ 
\sigma_u=\sqrt{1/M\sum_{a=1..M}(x_u^{(a)}-\bar{x}_u)^2} \nonumber
\end{align}

\subsection{Reduced Linear Response Modeling of the Controllable Resources}
\label{subsec:response_modeling}

The main idea of the proposed complexity reduction approach is to parameterize response of the controlled resources, $x_c$, which are primarily  generators, $Pg$, to the fluctuating component of the uncertainty characteristic, $w\doteq x_u-\bar{x}_u$. We will be using $x_c$ to represent all controlled variables and $Pg$ to represent a subset of control degrees of freedom associated with generator dispatch of active power. We consider the following (three) response policies/parametrizations, similar to the policies in \cite{roald2016integrated}:
\begin{itemize}
\item[a)] Response with a fixed (predefined) participation factor $\alpha$ of the controlled generators 
\begin{align}
& Pg_i(w)=\overline{Pg}_i+\alpha_i\cdot\Omega & i=PV,      \nonumber\\
& Pg_i(w)=\overline{Pg}_i+\alpha_i\cdot\Omega+\delta p(w) & i=\theta V,  \nonumber\\
& \Omega=\sum_{i \in N} w_i,\ 
\sum_{i \in N}\alpha_i=1.  \nonumber
\end{align}
Here $PV$ and $\theta V$ represent, respectively, set of buses where active power+voltage, and phase+voltage are kept fixed (the latter applies to the system's slack bus). 

\item[b)] A version of case a) where the linear response vector, $\alpha$, is not fixed, but is treated as an optimization variable. 

\item[c)] In the last, more general version, the linear response is parameterized by a matrix, $\alpha$, such that each generator $P_g$ responds separately to each component of $w$,
\begin{align}
& Pg_i(w)=\overline{Pg}_i+\sum_{j \in N} \alpha_{ij}\cdot\omega_j & i=PV     \nonumber\\
& Pg_i(w)=\overline{Pg}_i+\sum_{j \in N} \alpha_{ij}\cdot\omega_j+\delta p(w) & i=\theta V   \nonumber\\
& \sum_{i \in N}\alpha_{ij}=1 & \forall j \in N \nonumber
\end{align}
Here, $\alpha$ is an optimization variable. 
The response is still an affine, but allows for more general response patterns.
\end{itemize}
We refer to the three policies in the following uniform form:
\begin{align}
& Pg(w,\alpha)=Resp^{k}(w,\alpha) & k=a,b,c, \forall w \nonumber
\end{align}

\subsection{Cloud-AC-OPF}

We are now in the position to formulate the reduced model:
\begin{eqnarray}
&&\hspace{-0.5cm}\mbox{Cloud-AC-OPF}(\bar{x}_u,\sigma_u)=\min_{\bar{x}_c,\alpha} \mathbb{E}_{w}[Cost(\bar{x}_c+Resp^{k}(w,\alpha))] \nonumber\\ \hspace{-0.5cm}
&&\mbox{s.t. } \mbox{AC-PF}(\bar{x}_c,\bar{x}_u)=0\nonumber \\ \hspace{-0.5cm} && \mbox{Prob}_{w}\left[\mbox{Constr}(\bar{x}_c+\phi_c+Resp^{k},\bar{x}_u+w)\geq 0\right]\leq\varepsilon \label{eq:cloud-OPF}
\end{eqnarray}
where the expectation and the probabilistic expressions are given with respect to the uncertain/uncontrolled variable $w$. The functions $\phi_c=\phi_c(x_c,\bar{x}_u,w)$ describe variations of the controlled variables (other than $Pg(w,\alpha)$) as a function of the realization of $w$.

It is important to emphasize that although the formulation \eqref{eq:cloud-OPF} is similar to the CC-AC-OPF in, e.g., \cite{13ROKA,14BCH}, the relation between the Cloud-AC-OPF and the CC-AC-OPF is only formal.  As explained above, Cloud-AC-OPF represents a reduced model, where the linear response coefficient(s) $\alpha$ in (\ref{eq:cloud-OPF}) can be understood as a (conservative) approximation of the generators ability to react to uncertainty. On the contrary, the corresponding linear coefficients in CC-AC-OPF, see e.g. \cite{13ROKA,14BCH}, represent the actual automatic generation response to short term fluctuations. 

In the following we test how the three tractable reduction schemes, parameterized by $Resp^{k}(w,\alpha)$ with $k=a,b,c$ in Cloud-AC-OPF, approximate the MS-AC-OPF.
If the historical data is provided in the form of normal distributions with given properties then Cloud-AC-OPF can be directly applied.

\subsection{Cloud-AC-OPF: Analytic Reformulation}
When the exogenously introduced cloud of samples, representing fluctuations of the uncontrolled sources, $w$, around the center of the cloud, $\bar{x}_u$, is sufficiently small (or just moderate in size relative to the mean - fluctuations are small), we can linearize the non-linear AC power flow equations and still expect a reasonably accurate representation of the response to the fluctuations $\phi_c$. 
To obtain a tractable deterministic reformulation of the chance constraint in (\ref{eq:cloud-OPF}), we use a moment-based reformulation dependent only on the mean and standard deviation $\bar{x}_u, \sigma_u$. In fact, this reduction (tracking only two first moments) provides probabilistic guarantees for a much wider range of distributions with finite mean and variance \cite{roald2015security}, and can more generally be understood as a robust optimization with feasibility guarantees for uncertainty realizations  \cite{ben2009robust}, which still requires fixing safety level parameters. 

With the assumptions of linearized AC power flow equations and a moment-based chance constraint reformulation, we arrive at the following version of (\ref{eq:cloud-OPF})
\begin{eqnarray}
&&\min_{\bar{x}_c,\alpha} \mathbb\mbox{Cost_1}(\bar{x}_c)+\mathbb{E}_w\left[Cost_2(w,\alpha)\right] \label{eq:cloud-OPF-simpler} \\ \hspace{-0.5cm}
&&\mbox{s.t. } \mbox{AC-PF}(\bar{x}_c,\bar{x}_u)=0\nonumber \\ \hspace{-0.5cm} && \mbox{Prob}_{w}\left[\mbox{Constr}(\bar{x}_c+G_w w+Resp^{k}(w,\alpha))\geq 0\right]\leq\varepsilon  \nonumber
\end{eqnarray}
where the adjustment of the controlled variables $\phi_c$ is defined using sensitivity matrix $G_w$, describing linear response of the controlled variables to variations in the exogenous/uncontrolled variables:
\begin{eqnarray}
\phi_c =\left.\frac{\partial x_c}{\partial w }\right|_{\begin{array}{l}x_c=\bar{x}_c\\ x_u=\bar{x}_u\end{array}}w \doteq  G_w(\bar{x}_c,\bar{x}_u,\alpha)w. \label{eq:lin} \nonumber
\end{eqnarray}
If AC-PF system equations are feasible in the center of the cloud (at $w=0$) then the sensitivity matrix and $\phi_c$ exists and differentiable.
Here in (\ref{eq:cloud-OPF-simpler}) explicit expression for $G_w$, as a function of $\bar{x}_c,\bar{x}_u,\alpha$, is skipped due to space limitations; the objective function is split in two parts, correspondent to mean and fluctuations, respectively.  

Following the approach of \cite{13ROKA,14BCH,18RA}, we are able to evaluate the expectation and the probabilities in (\ref{eq:cloud-OPF-simpler}) analytically. Moreover, the analytic evaluation returns explicit dependencies on $\bar{x}_c$ and $\alpha$,  therefore stating the Cloud-AC-OPF (\ref{eq:cloud-OPF-simpler}) as the following tractable deterministic optimization formulation:
\begin{eqnarray}
&&\min_{\bar{x}_c,\alpha} \mathbb\mbox{Cost_1}(\bar{x}_c)+\mathbb{E}_w\left[Cost_2(w,\alpha)\right] \label{eq:cloud-OPF-1}\\ \hspace{-0.5cm}
&&\mbox{s.t. } \mbox{AC-PF}(\bar{x}_c,\bar{x}_u)=0\nonumber \\ 
\hspace{-0.5cm} && \mbox{Constr}(\bar{x}_c, \bar{x}_u)\leq -\lambda(\bar{x}_c, \bar{x}_u,\alpha, \Sigma_w), \nonumber 
\end{eqnarray}
where the dependence of the correction to the cost on $\alpha, \bar{x}_c$ and samples is detailed below. Uncertainty margins, $\lambda(\bar{x}_c, \bar{x}_u,\alpha, \Sigma_w)$, are computed for each type of variables ($\gamma$) and each type of bus/line ($\mu$: PQ,PV,$\theta V$,line) as:
\begin{align}
& \lambda_{\gamma:\mu} =0, {\scriptstyle V:\theta V,V:PV,P:PQ,Q:PQ} \nonumber \\
& \lambda_{ P:PV/\theta V}^i = \Phi^{-1}(1-\epsilon_{\gamma}) \times ||(G^{\gamma:\mu}_{w(i,:)}+\alpha_{(i,:)}) \Sigma_w^{1/2}||_2 \nonumber \\
& \lambda_{\gamma:\mu}^i = \Phi^{-1}(1-\epsilon_{\gamma}) \times ||(G^{\gamma:\mu}_{w(i,:)}) \Sigma_w^{1/2}||_2, {\scriptstyle V:PQ,Q:PV/\theta V,F:line}  \nonumber
\end{align}
where $\Phi^{-1}$ stands for the inverse cumulative distribution function of the standard normal distribution. $\Sigma_w$ is the calculated covariance matrix for uncontrolled sources for the cloud. F is squared apparent power at from/to side of a line. The resulting deterministic optimization (\ref{eq:cloud-OPF-1}) over $\bar{x}_c$ and $\alpha$ does not have any nice structure (e.g. it is not convex). 

\subsection{Analytic Averaging of Cost over the Uncertainty Set}
In the case of quadratic dependence of the objective on the generation dispatch,  the fluctuating part of the cost becomes 
\begin{align}
& {\scriptstyle i-PV:} \mathbb{E}_w\left[Cost_2^i(w,\alpha)\right]=a\sum_{j=1..N}\alpha_{ij}^2 v(w_j) \nonumber \\
& {\scriptstyle i-\theta V:} \mathbb{E}_w\left[Cost_2^i(w,\alpha)\right]=a\sum_{j=1..N}(\alpha_{ij}+(G^{P:\theta V}_w)_{1j})^2 v(w_j) \nonumber
\end{align}
where formal expectation over $\omega$ is stated in terms of variances of the uncertainty, $v(w_j)$, at the $j$-th uncertainty site, evaluated directly from (available) samples; and $G^{p:\theta V}_w$ is 1-row submatrix of $G_w$ corresponding to sensitivity of active power at slack bus to $w$.

\subsection{Implementation and Solution Approach}

We solve (\ref{eq:cloud-OPF-1}) via
iterative algorithm  implemented in Julia using JUMP \cite{JUMP}, thus taking advantage of the modularity and automatic differentiation features of the software. The idea of the algorithm is to specify (current) $x\doteq (\bar{x}_c, \alpha)$ at each iteration step. Then the sensitivity matrices $G_w$ are evaluated at $\bar{x}_c$ with analytical dependence only on $\alpha$ then provided as an input to the optimization model. This dependence also applies to the uncertainty margins, $\lambda$. 
Schematic description of the algorithm is as follows:
\begin{itemize}
\item[1.] Initialization. Set uncertainty margins $\lambda_P^0=\lambda_Q^0=\lambda_V^0=\lambda_F^0=0$. Solve classical AC-OPF, and set its argmin as $\bar{x}_{c}^0$. Set iteration number to $k=1$.

\item[2.] Evaluate starting point for optimization variables. $\alpha^{start}$ - defined as the equal participation of each generator evaluated for each component of the uncertainty vector. $\lambda_P^{start},\lambda_Q^{start},\lambda_V^{start},\lambda_F^{start}$ are computed using sensitivity matrices at $\bar{x}_{c}^{k-1}$ and $\alpha^{start}$. $\bar{x}_{c}^{start}$ and auxiliary variables - take $k-1$ solution. 

\item[3.] Define non-linear optimization model (\ref{eq:cloud-OPF-1}). Model variables are $\bar{x}_c=$ $\bar{V}_c,\bar{\theta}_c,\bar{Pg}_c,\bar{Qg}_c;\alpha;\lambda_P^{var};\lambda_Q^{var};\lambda_V^{var};\lambda_F^{var}$ and $\overline{p_{fr}}_c;\overline{p_{to}}_c;\overline{q_{fr}}_c;\overline{q_{to}}_c$ - the auxiliary variables. Define constraints according to (\ref{eq:cloud-OPF-1}). Effectively model is similar to AC-OPF but with corrected by the uncertainty margins constraints, additional variables and modified objective - averaged over the cloud.

\item[4.] Solve the model. Update $k=k+1$. Go to step 2.
\end{itemize}

\section{The Case Study Analysis}

We test our approach on the IEEE 30 bus system available within the  Matpower package \cite{MATPOWER}. We set the base configuration, $\bar{x}_u$, of the uncontrolled degrees of freedom and draw samples for active power consumption/production at the uncertain node from a Gaussian distribution with a prescribed covariances (standard deviations) around the base case. Components of the vector of standard deviations is defined according to a prescribed ratio of the base case loading (for example 5$\%$ of initial active demand). Reactive power is set constant (base case). Uncertainty set size of 1000 samples was defined experimentally to provide stable MS-AC-OPF cost for different sample sets at given base case and level of uncertainty. We solve the MS-AC-OPF directly for all the samples, therefore setting up ground-truth standards for  the following comparisons. Then we solve different version of the Cloud-AC-OPF model and compare the solutions with the ground truth set by the MS-AC-OPF solution. 

According to the introduced response model, voltage set-points on generators are fixed for different samples. Fixed reactive power does not reduce generality of the approach, if Q is uncertain it would be taken into account by similar sensitivity/uncertainty margin means.




\subsection{MS-AC-OPF and Three Flavours of Cloud-AC-OPF}

We experiment with the three flavours of the Cloud-AC-OPF model described in Section \ref{subsec:response_modeling}. In the following, we refer to the MS-AC-OPF with 1000 samples as MS-AC-OPF, and models (a-c) introduced in Section \ref{subsec:response_modeling} as Cloud-AC-OPF-k (k=a,b,c).

\subsection{Scalability Analysis}

Computations are done on a 3.3 GHz core i7 laptop CPU. The results are summarized in the Table \ref{tab1}:

\begin{table}[h!]
\caption{Computational time comparison for the introduced models} 
\centering 
 \begin{tabular}{|c| c| c| c|} 
 \hline
 Model & Description & $\#$ Iterations & Time (sec) \\ [0.5ex]
  \hline
 MS-AC-OPF & 1000 samples & 1 & 70 \\ 
 \hline
 Cloud-AC-OPF-a & given $\alpha$ & 5 & 2.8 \\ 
 \hline
 Cloud-AC-OPF-b & $\alpha$-vector & 5 & 41.6 \\
 \hline
 Cloud-AC-OPF-c & $\alpha$-matrix & 5 & 41.4 \\ 
 \hline
\end{tabular}
\label{tab1} 
\end{table}

All the Cloud-AC-OPF cases are solved faster than the bare MS-AC-OPF. We expect to see this acceleration effect to be even more pronounced in the case of the aforementioned two-stage planning models. Also MS-AC-OPF computational time grows linearly with number of samples while Cloud-AC-OPF time depends only on the system size.

\section{Reduction Analysis}

In this section we discuss quality of the Cloud-AC-OPF solutions as they are compared to the benchmark provided by the MS-AC-OPF.  We choose the optimal objective function and the optimal active generation dispatch as benchmarks for comparison. 

The following two operational regimes are considered:
\begin{enumerate}
\item Low loading - in this case all line, voltage and power injection conditions are safely within the feasible operational limits (not saturated). 
\item High loading - this is a heavily congested case with a number of constraints being either saturated or close to be saturated.  
\end{enumerate}

\subsection{Comparison by the Value of Optimal Objective}

\begin{figure}[h!]
\centering
\includegraphics[width=0.45\textwidth]{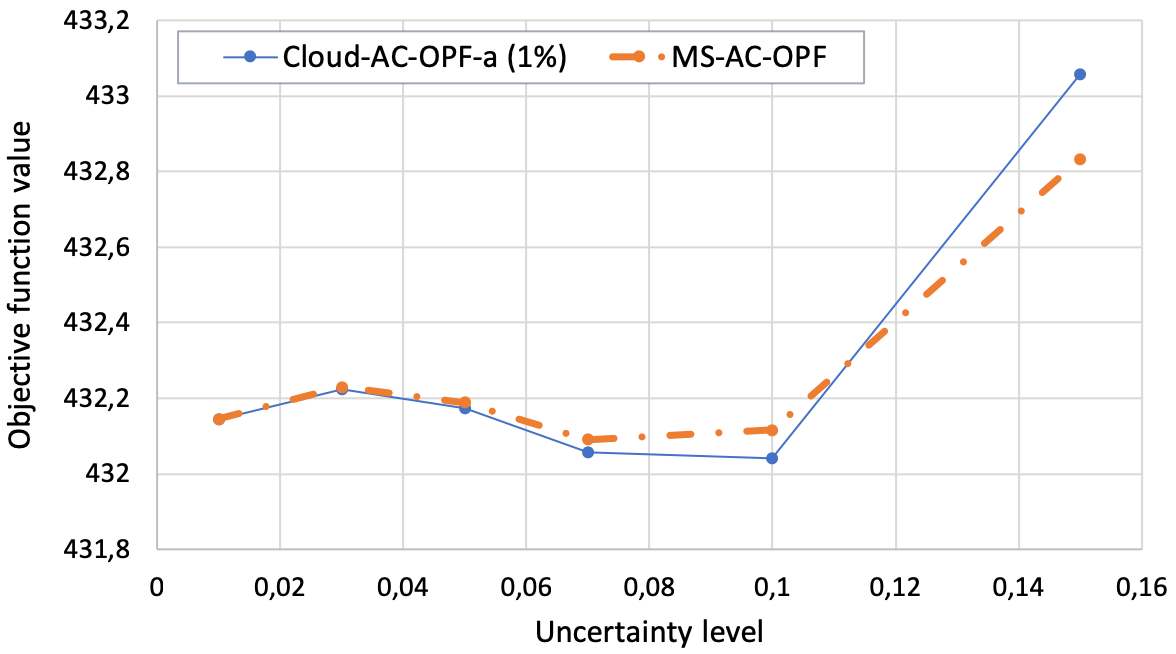}
\caption{Obj. function comparison for low loading regime (MS-AC-OPF vs Cloud-AC-OPF-a). Probability of the chance constraints violation is $\epsilon=1\%$}
\label{low_loading_cost}
\end{figure}

Consider, first, the low-loading regime - the base case loading is set to $0.8$ and the probability of a constraint violation is set to 1$\%$ (for all the constraints). The results are shown in Fig.~\ref{low_loading_cost}. In this case  even the least accurate Cloud-AC-OPF-a shows a satisfactory performance. We observe a good approximation quality up to 10$\%$ of the uncertainty level, when assessed by the objective function value. 

\begin{figure}[h!]
\centering
\includegraphics[width=0.45\textwidth]{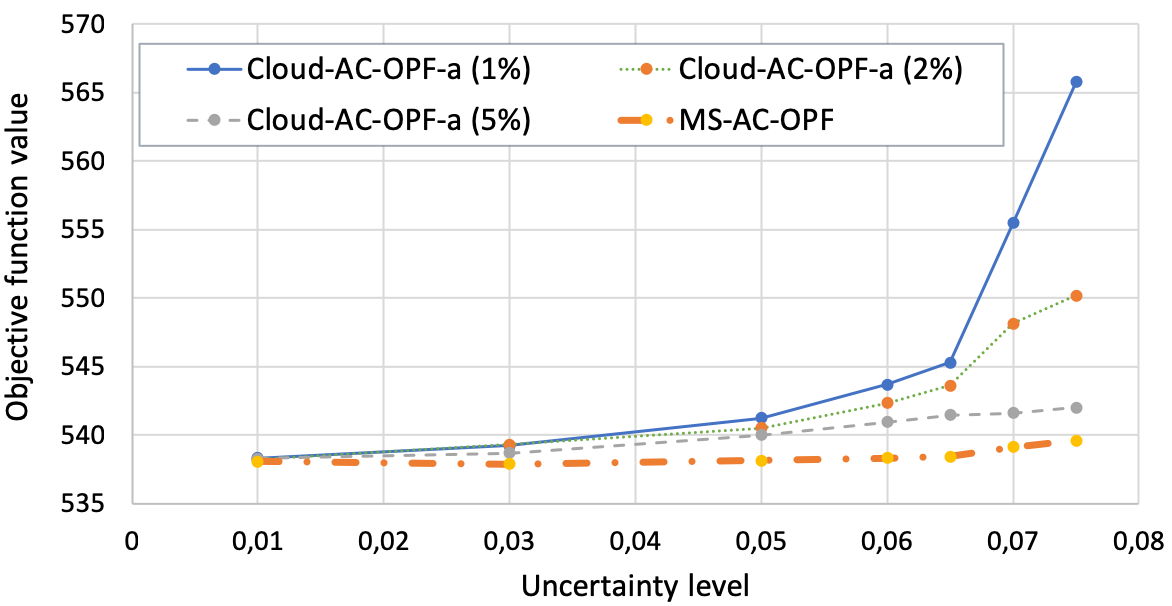}
\caption{Obj. function comparison in the regime of high load (MS-AC-OPF vs Cloud-AC-OPF-a) shown as a function of the uncertainty level and $\epsilon$.}
\label{high_loading_eps}
\end{figure}

In the high loading regime, when the base case loading is set to $0.95$ (the loading level $1.03$ would be already AC-OPF infeasible), a gap in performance between the MS-AC-OPF and the Cloud-AC-OPF appears.  The gap increases with the uncertainty (variance of fluctuations). Fig.~\ref{high_loading_eps} illustrates this effect for the Cloud-AC-PF-a evaluated with the probability of constraints violation from 1$\%$ to 5$\%$. 

It is also observed that the Cloud-AC-OPF models highly sensitive to variations of $\epsilon$ (probability of the chance constraints violation) in the high loading and high uncertainty regimes. This would be important when we will work with historical data. Uncertainty distribution is not necessary Gaussian in that case. It would be approximated by Gaussian with computed mean and variance. And $\epsilon$ should be carefully chosen (e.g. by doing out of sample analysis of the solution).

\begin{figure}[h!]
\centering
\includegraphics[width=0.45\textwidth]{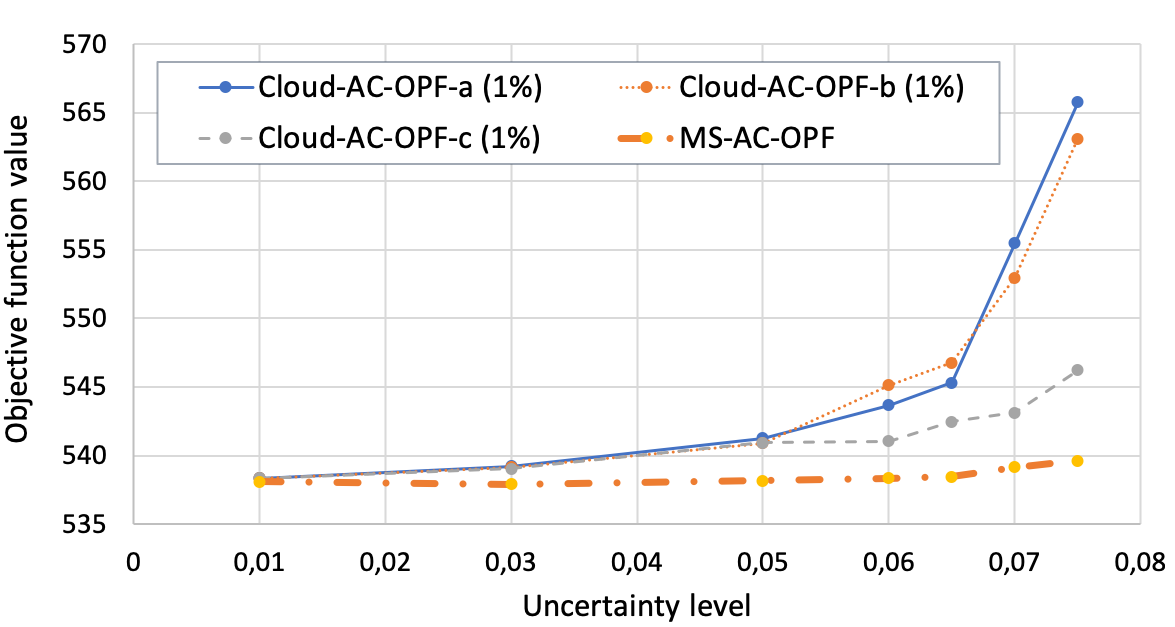}
\caption{Obj. function comparison at high loading regime for different versions of the Cloud-AC-OPF-k. $\epsilon=1\%$ for each model.}
\label{high_loading_alg}
\end{figure}

Fig.~\ref{high_loading_alg} compares objective function value of the Cloud-AC-OPF-k with the MS-AC-OPF depending on the generators response parametrization model. Probability of constraints violations is set to 1$\%$. It can be observed that Cloud-AC-OPF-b/c demonstrate better performance than the simplest model and that the most sophisticated matrix version Cloud-AC-OPF-c is advantageous in the high loading and high uncertainty regimes.  

Configurational (generation dispatch) comparison is performed in the following subsection.

\subsection{Comparison of the Cloud-AC-OPF and MS-AC-OPF Optimal Dispatches}

Comparison of optimal dispatches in the original MS-AC-OPF model and the reduced Cloud-AC-OPF model constitutes a much richer test (than based on the optimal cost) thus setting better criteria for the assessment of model reduction. The configuration (optimal dispatches based) analysis of the Cloud-AC-OPF shows significant dependence on the loading regime and also strong sensitivity to selection of the type of response in the model (a-c). 

\begin{figure}[h!]
\centering
\includegraphics[width=0.5\textwidth]{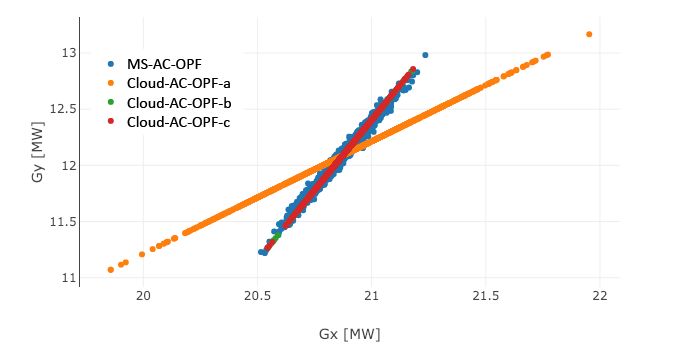}
\caption{Illustration of the Cloud-AC-OPF solution quality in terms of the generation response to perturbations in the regime of low loading.}
\label{low_loading_fixed}
\end{figure}

Fig.~\ref{low_loading_fixed} shows details of the comparison of Cloud-AC-OPF-k ($k=a-c$) with the MS-AC-OPF in the low-loading regime. 

The figure visualizes response of a representative pair of generators. Output from MS-AC-OPF model is shown in terms of samples (blue dots). Affine response of different versions of the Cloud-AC-OPF-k is shown for the comparison. In Cloud-AC-OPF-a (orange) all generators participate equally which is not relevant to actual shape of the MS-AC-OPF state space. In Cloud-AC-OPF-b/c (green and purple) the response is optimized. Original state space is close to $(k=b)$ model in the low loading regime (gens effectively respond to total power mismatch). Because of that both Cloud-AC-OPF-b/c demonstrate good performance and matrix response model basically finds the same as vector response model.   

\begin{figure}[h!]
\centering
\includegraphics[width=0.5\textwidth]{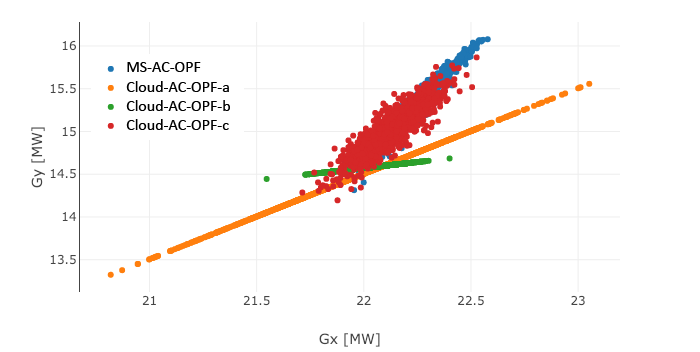}
\caption{Illustration of the Cloud-AC-OPF solution quality in terms of the generation response to perturbations in the regime of high loading.}
\label{high_loading_fixed}
\end{figure}

In the regime of a high load, illustrated in Fig.~\ref{high_loading_fixed}, some constraints become active. This results in the fact that MS-AC-OPF state space of optimal configurations/dispatches is more complicated.  We observe that in this case the more constrained models Cloud-AC-OPF-a/b (green and orange) fails to represent MC-AC-OPF. However, one also observes that solution quality of the more advanced model reduction scheme, represented by Cloud-AC-OPF-c (purple), is still satisfactory. (Here, in the case of Cloud-AC-OPF-c, probability of constraints violation over lines is increased to 10$\%$, while probability of all other violations is kept to the 5$\%$). 

\section{Conclusion}

The material presented in this manuscript demonstrates that the cloud-based approach has a strong potential in the model reduction applications. The significant reduction in the computational effort is achieved due to the fact that the Cloud-AC-OPF modeling, implemented via generalized chance-constraints, allows to represent infinitely many realistic configurations with a small number of decision variables.
Resulting solutions demonstrate good approximation quality,  when compared with the ground truth, set by the MS-AC-OPF, in terms of the optimal objective. The best approximation quality is achieved when the reduced model is represented by a matrix response,  where each generator responds to each variable source (i.e., each network node with significant variability) independently. 
These first results are encouraging as they open a multitude of opportunities to handle more complex and challenging settings, such as two-level and multi-stage planning problems of the type discussed in \cite{93GCCP,2015MW,17MPSS}.

We conclude with a summary of technical contributions reported in the manuscript: 
\begin{enumerate}
\item We develop the Cloud-AC-OPF approach, which applied the CC-OPF methodology \cite{13ROKA,14BCH,18RA} to model reduction of the computationally prohibitive MS-AC-OPF. 

\item The approach takes advantage of an affine parametrization of the decision variables. 

\item We show practicality of the approach on the examples of moderate size (IEEE 30 bus model) where validation against the MS-AC-OPF, providing the ground truth, is still feasible.
\end{enumerate}

We plan to extend this work in the future to multi-cloud situations representing richer and more realistic historical data (i.e. uncertainty which cannot be represented by a single quasi-Gaussian cloud) examples, relevant to resolving most difficult two-level and multi-time frame expansion planning problems.

\bibliographystyle{IEEEtran}
\bibliography{cloud.bib}

\end{document}